\documentclass[letterpaper,twoside,reqno]{amsart}

\usepackage{epsf,graphicx} 

\usepackage{latexsym,amsfonts}

\newcommand{\conv}{\operatorname{conv}}

\newcommand{\R}{{\mathbb R}}
\newcommand{\A}{{\mbox{$\mathcal A$}}}
\newcommand{\HH}{{\mbox{$\mathcal H$}}}

\newtheorem{theorem}{Theorem}[section]
\newtheorem{proposition}[theorem]{Proposition}
\newtheorem{lemma}[theorem]{Lemma}
\newtheorem{corollary}[theorem]{Corollary}

\newtheorem{definition}[theorem]{Definition}

\newtheorem{remark}[theorem]{Remark}
\newtheorem{remarks}[theorem]{Remarks}
\newtheorem{example}[theorem]{Example}

\title{The polytope of non-crossing graphs on a planar point set}

\author{David Orden and Francisco Santos}

\address{Dpto. de Matem\'{a}ticas, Estad\'{\i}stica y Computaci\'on.
Universidad de Cantabria. 39005 Santander, SPAIN}

\email{ \{ordend,santos\}@matesco.unican.es}

\date{May, 2003}

\thanks{Research partially supported by project BMF2001-1153 of
Spanish Direcci\'on
General de Investigaci\'on Cient\'{\i}fica.}

\begin{document}
\thispagestyle{empty}


\begin{abstract}
For any finite set $\A$ of $n$ points in $\R^2$,
we define a $(3n-3)$-dimensional
simple polyhedron whose face poset is isomorphic to the poset of
``non-crossing marked  graphs''  with vertex set $\A$,
where a marked graph is defined as a geometric graph together with a
subset of its vertices.
The poset of non-crossing graphs on $\A$ appears as the complement of
the star of a face in that polyhedron.

The polyhedron has a unique maximal bounded face, of dimension
$2n_i +n -3$ where $n_i$ is the number of points of $\A$ in the interior of
$\conv(\A)$.
The vertices of this  polytope are all the
pseudo-triangulations of $\A$, and the edges are flips of two types:
the traditional
diagonal flips (in pseudo-triangulations) and the removal or
insertion of a single edge.

As a by-product of our construction we prove that all pseudo-triangulations are
infinitesimally rigid graphs.

\end{abstract}

\maketitle

\section{Introduction and statement of results}

\subsection*{The polyhedron of non-crossing graphs}
The set of (straight-line, or geometric) non-crossing graphs with a given
set of vertices in the plane is of interest in Computational
Geometry, Geometric
Combinatorics, and related areas. In particular, much effort has been
directed towards
enumeration, counting and optimization on the set of maximal such
graphs, that is to say,  triangulations.


The poset structure of the set of non-crossing graphs, however,   is only well understood
if the points are in convex position. In this case the non-crossing graphs
containing all the hull edges are the same as the polygonal
subdivisions of the convex
$n$-gon and, as is well-known, they form the face poset of the
\emph{$(n-3)$-associahedron}.

In this paper we  construct from $\A$ a
polytope  whose face poset contains the poset of non-crossing graphs
on $\A$ embedded in a very nice way. In the following statement, and in
the rest of the paper, we  distinguish between vertices of $\conv(\A)$ 
(\emph{extremal} points) and other points in 
the boundary of
$\conv(\A)$, which we call \emph{semi-interior} points. 

\begin{theorem}
Let $\A$ be a finite set of $n$ points in the plane, not all 
contained in a line.
Let $n_i$, $n_s$ and $n_v$ be the number of interior, semi-interior and
extremal points of $\A$, respectively.

There is a simple polytope $Y_f(\A)$ of dimension
$2n_i + n - 3$, and a face $F$ of $Y_f(\A)$ (of dimension $2n_i + n_v - 3$)
such that the
complement of the star of $F$ in the face-poset of $Y_f(\A)$
equals the
poset of non-crossing graphs on $\A$ that use all the convex hull edges.
\label{thm:posetintro}
\end{theorem}

This statement deserves some words of explanation:
\begin{itemize}

\item[--] Our equality of posets reverses  inclusions.
Maximal non-crossing graphs (triangulations of $\A$) correspond to
minimal faces (vertices of $Y_f(\A)$).

\item[--] We remind the reader that the \emph{star} of a
face $F$  in a polytope is the subposet of facets that contain $F$ and all their faces.
In the complement of the star of $F$ we must consider the polytope
$Y_f(\A)$ itself as an element, which corresponds to the graph with no interior edges.

\item[--] If $\A$ is in convex position, then $Y_f(\A)$ is 
the associahedron. The face $F$ is the whole polytope, whose star
we must interpret as being empty.

\item[--] Our results are valid for point sets in non-general 
position. Our definition of non-crossing in this case
is that if $q$ is between $p$ and $r$, then the edge $pr$ cannot 
appear in a non-crossing
graph, regardless of whether $q$ is incident to an edge in the 
graph. This definition
has the slight drawback that a graph which is non-crossing in a point 
set $\A$ may become
crossing in $\A\cup\{p\}$, but it makes maximal crossing-free graphs
coincide with the triangulations of $\A$ (with the convention, 
standard in Computational Geometry,
that triangulations of $\A$ are required to use all the points of 
$\A$ as vertices). To be consistent with this,
a convex hull edge is an edge between
two consecutive extremal or semi-interior
points of $\A$ in the boundary of $\conv(\A)$.

\item[--] Since convex hull edges are irrelevant to crossingness,
the poset of {\em all} non-crossing graphs on $\A$
is the direct product of the poset in the statement
and a Boolean poset of rank $n_v + n_s$.
\end{itemize}

We give a fully explicit facet description of $Y_f(\A)$. It lives in
$\R^{3n}$ and is defined by the 3 linear equalities (\ref{eq:S}) and
the ${n\choose 2}+n$ linear inequalities (\ref{eq:vf}) and 
(\ref{eq:tf}) of Section
\ref{section:polyhedron}, with some of them turned into equalities. 
(With one exception:
for technical reasons,
if $\A$ contains three collinear boundary points we need to add extra
points to its exterior and obtain $Y_f(\A)$ as a face of the polytope
$Y_f(\A')$ of the extended point set $\A'$).

The $f_{ij}$'s in equations
(\ref{eq:vf}) and (\ref{eq:tf})  and in the  notation $Y_f(\A)$ 
denote a vector in
$\R^{n+1\choose 2}$. Our construction starts with a linear cone 
$\overline{Y_0}(\A)$
(Definition \ref{defi:cone}) whose facets are then translated
using the entries of $f$ to produce a polyhedron $\overline{Y_f}(\A)$, of
which the polytope $Y_f(\A)$ is the unique maximal bounded face.
Our proof goes  by analyzing the
necessary and sufficient conditions for $f$ to produce a polytope
with the desired
properties and then proving the existence of valid choices of $f$. In
particular, Theorem \ref{thm:exist_valid_f} shows one valid choice.
This is essentially the same approach used in
\cite{Ro-Sa-St} for the polytope of pointed non-crossing graphs 
constructed there, which
is actually the face $F$ of $Y_f(\A)$ referred too in the statement.
Our results generalize the construction in \cite{Ro-Sa-St}
in two directions: on the one hand they refer to all non-crossing graphs
(not only pointed ones) and on the other hand we show how to deal
with non-general position.



Mapping the instances of a combinatorial object to the vertices of a certain polytope
is useful for optimization and enumeration purposes. In the case of triangulations, two such 
polytopes have been used in the past; the  ``secondary polytope'' \cite{Bi-Fi-St}
and the ``universal polytope'' \cite{Lo-Ho-Sa-St}
of the point set. Our construction adds to these two, but it has one 
advantage: We have an explicit facet description of the polytope. In the secondary polytope, facets
correspond to the coarse polygonal subdivisions of $\A$,
which have no easy characterization. In the universal polytope, the
facet description in \cite{Lo-Ho-Sa-St} gives only a linear programming relaxation
of the polytope, which makes
integer programming be needed in order to optimize linear functionals on it.

That the poset we are interested in is not the whole poset of faces of the polytope
$Y_f(\A)$
may seem a serious drawback for using it as a tool for enumeration  of  all the triangulations of
a planar point set. But the subposet we are interested in is not just a subposet.
Being the complement of  the  star of a face $F$ has 
theoretical and practical implications. On the one hand, it
implies that the poset is a shellable  ball of dimension  $2n_i+n-4$,
since there is a shelling order ending in the facets that
contain $F$. On the other hand, the part of the boundary of $Y_f(\A)$ 
that we are
interested in becomes the (strict) lower envelope of a convex polyhedron
via any projective transformation that  sends a supporting hyperplane
of $F$ to the infinity.


\subsection*{The polytope of pseudo-triangulations}
Actually, the set of {\em all} the vertices of the polytope $Y_f(\A)$ is interesting on its
own merits:

\begin{theorem}
\label{thm:polytopeintro} The vertex set of the polytope $Y_f(\A)$ of Theorem
\ref{thm:posetintro} is in bijection to  the set of all
pseudo-triangulations of $\A$.
The 1-skeleton of $Y_f(\A)$ is the graph of flips between them.
\end{theorem}

\begin{figure}[htb]
\begin{center} {\ }
\epsfxsize=3in\epsfbox{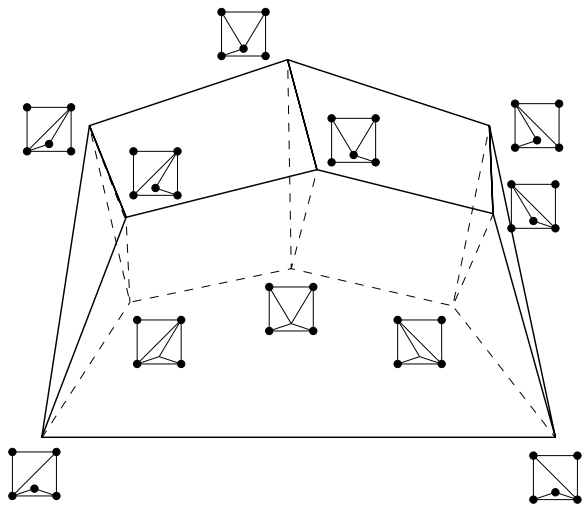} \\
\caption{The pseudo-triangulations of a point set 
form a polytope}\label{fig:graph_sq}
\end{center}
\end{figure}

A pseudo-triangle is a perhaps non-convex polygon with three corners (non-reflex vertices).
A pseudo-triangulation of $\A$ is a subdivision of $\conv(\A)$ into pseudo-triangles
using $\A$ as vertex set. 
As an example, Figure \ref{fig:graph_sq} shows the pseudo-triangulations of a set of 
four extremal points
and one interior point,
with the graph of flips between them embedded as (a Schl\"egel diagram of) a 
4-dimensional polytope.  The face $F$ of Theorem
\ref{thm:posetintro} is actually the facet on which the Schl\"egel diagram has been 
performed. Why the interior point is ``unmarked'' in the three triangulations
of the point set and ``marked'' in the other 8 pseudo-triangulations will be explained in
Section \ref{section:pseudot}.

Pseudo-triangulations, first introduced by Pocchiola
and Vegter around 1995 (see \cite{PV}), have by now been used in many
Computational
Geometry applications, among them visibility
\cite{PV-theory, PV2, PV,Toth-Speck}, ray shooting \cite{gt}, and kinetic data
structures \cite{guibas, KSS}. Streinu \cite{Streinu} introduced the
\emph{minimum} or \emph{pointed pseudo-triangulations} (p.p.t.'s, for short), 
proved certain surprising relations of them  to structural
rigidity of non-crossing graphs, and used  them to prove the
Carpenter's Rule
Theorem (the first proof of which was given shortly before by  Connelly et al
\cite{Co-De-Ro}). P.p.t.'s can also be described as the 
maximal non-crossing  and pointed graphs. See all the relevant definitions 
in Section \ref{section:pseudot}.

The paper  \cite{Ro-Sa-St},  using these rigid theoretic properties of p.p.t.'s, 
constructed a polytope $X_f(\A)$ whose vertex set is precisely the set of p.p.t.'s
of a given point set. Our construction is heavily based on
it (but this paper is mostly self-contained). Our method not only
extends that construction to cover all pseudo-triangulations, it also shows that
the method can easily be adapted to point sets in non-general position,
with a suitable definition of pseudo-triangulation for such point
sets  (Definition \ref{defi:pseudot-sp}).

The flips between pseudo-triangulations that appear in our
Theorem \ref{thm:polytopeintro} are
introduced in Section
\ref{section:pseudot} (see Definition \ref{defi:ptflips}) for point sets
in general position and in Section \ref{section:special}
(see Definition \ref{defi:ptflip-sp} and Figures
\ref{fig:5pts} and \ref{fig:6pts})
for point sets with collinearities. The definition  is new (to the
best of our knowledge) but in the case of general position it
has independently been considered in
\cite{AABK}, where flips between pseudo-triangulations are related to
geometric flips between polyhedral terrains. It is proved there that the
graph of flips of a point set $\A$ (that is, the 1-skeleton of our
polytope) has diameter bounded by $O(n\log(n))$. The same bound for 
the graph of pointed pseudo-triangulations (the 1-skeleton of the
face $F$) has been obtained in \cite{Besp}.

Our flips restrict to the ones in \cite{Streinu} and \cite{Ro-Sa-St}
when  the two
pseudo-triangulations involved are pointed, and they are
also related to the flips of Pocchiola and Vegter \cite{PV} as
follows: Pocchiola and
Vegter were interested in pseudo-triangulations of a set ${\mathcal
O}:=\{o_1,\dots,o_n\}$ of  convex bodies, and they defined a graph of
flips between them. That graph is regular of degree
$3n-3$. Pocchiola (personal communication) has shown  that
taking each $o_i$ to be a sufficiently small convex body around  each 
point, our
graph is obtained
from the one in  \cite{PV} by contraction of certain edges.

\subsection*{Constrained pseudo-triangulations}

It is sometimes convenient to study pseudo-triangulations or non-crossing graphs on $\A$
that contain a specified subset of edges. For example, the
pseudo-triangulations of a non-convex polygon $P$, are in bijection with 
pseudo-triangulations of the point set
$\A=\operatorname{vertices}(P)$ which use all the boundary
edges of $P$ and any (arbitrarily chosen) triangulation of $\conv(P)\backslash P$.
We can even allow for additional interior points of $P$ to be used as vertices
(a situation called an ``augmented polygon'' in \cite{AABK}).

We want to emphasize that our results apply to
these and other cases: since each facet of $Y_f(\A)$ corresponds
to the presence of a certain edge, pseudo-triangulations or non-crossing graphs
that contain a specific non-crossing graph $G$ correspond to faces or vertices of 
$Y_f(\A)$ contained in the face that corresponds to $G$.

\begin{corollary}
\label{coro:polygon}
Let $G$ be a non-crossing graph on $\A$ with $k$ non-boundary edges 
(edges not in the boundary of $\conv(\A)$). Then, the set of pseudo-triangulations of $\A$
that contain $G$ as a subgraph is the vertex set of a simple polytope of dimension 
$2n_i+n-3-k$, whose 1-skeleton is the graph of flips between them.

In particular, the pseudo-triangulations of a non-convex polygon $P$ with $n$ vertices 
are the vertex set of a polytope of dimension $n-3$.
\end{corollary}

\begin{figure}[htb]
\begin{center} {\ }
\epsfxsize=3in\epsfbox{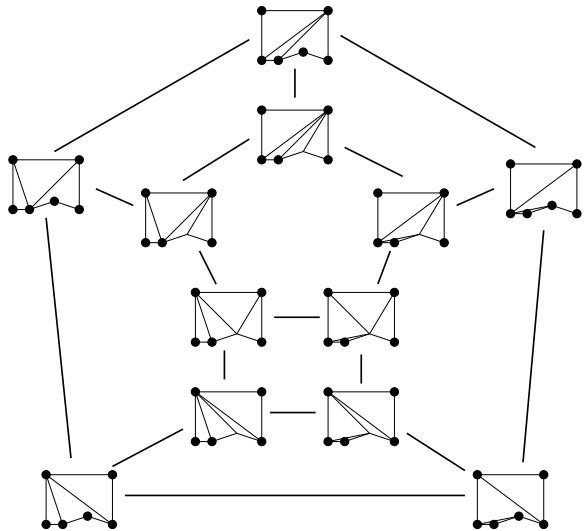} \\
\caption{The pseudo-triangulations of a (perhaps non-convex) polygon
form a polytope too}\label{fig:graph_hex}
\end{center}
\end{figure}

Figure \ref{fig:graph_hex} shows the (3-dimensional) polytope of pseudo-triangulations of a
certain non-convex hexagon.

\subsection*{Pseudo-triangulations are rigid}

Our construction has also rigid-theoretic consequences. A generically 
rigid graph (for dimension 2) 
is a graph which becomes rigid in almost all its straight-line embeddings in the plane. 
Generically rigid graphs need at least $2n-3$ edges, because that is the number of degrees of freedom
of $n$ points in the plane (after neglecting rigid motions). Generically rigid graphs with exactly $2n-3$ 
edges are called \emph{isostatic} and they admit the following characterization, due to Laman
(see, for example, \cite{gss}): 
they are
the graphs with $2n-3$ edges and with the property that any subset of $k\le n-2$ vertices is incident to at least 
$2k$ edges. Using this characterization, 
Ileana Streinu \cite{Streinu} proved that every pointed pseudo-triangulations is an isostatic graph. We have the following generalization:

\begin{theorem}
\label{thm:rigidintro} Let $T$ be a pseudo-triangulation of a planar point set $\A$ 
in general position. Let $G$ be its underlying graph.
Then:
\begin{enumerate}
\item $G$ is infinitesimally rigid, hence rigid and generically rigid.
\item There are at least $2k+3l$ edges of $T$ incident to any subset of $k$ pointed plus $l$ non-pointed
vertices of $T$ (assuming $k+l\le n-2$).
\end{enumerate}
\end{theorem}

In the statement, ``general position'' can actually be weakened to ``no three boundary points of
$\A$ are collinear.  In the presence of boundary collinearities, non-rigid pseudo-triangulations (for our definition) exist. For example, only six of the fourteen 
pseudo-triangulations of the point set of Figure \ref{fig:6pts} are rigid.

If we recall that a pseudo-triangulation with $k$ non-pointed vertices has exactly $2n-3+k$ edges
(see Proposition \ref{prop:pseudot}), Theorem \ref{thm:rigidintro}
 has the consequence that the space of self-stresses
on a pseudo-triangulation has exactly dimension $k$. This fact follows also from the results of \cite{AABK}.

\medskip




\subsection*{Two open questions}

\begin{itemize}

\item Can every planar and generically rigid graph be embedded as a pseudo-triangulation? 
That the answer is positive should be considered a conjecture of the set of authors of  \cite{barbados}, a paper that actually solves the minimal case: planar Laman graphs can be embedded as pointed pseudo-triangulations. The maximal case
(combinatorial triangulations can be drawn with convex faces) is also solved, by Tutte's theorem 
\cite{Tutte}.
But the intermediate  cases remain open. 

\item Are non-crossing graphs on $\A$ the face poset of a polyhedron? A naive answer would be
 that such a polyhedron is obtained by just deleting from the facet definition of 
$Y_f(\A)$  those facets containing $F$. But this does not
work even  in the point set with a single point in general position in the interior of a quadrilateral
(the one in Figure \ref{fig:graph_sq}). Using the equations that define the polytope it
can be checked that the removal of the facet $F$ in this example
gives a polyhedron with two extra vertices, not present in $Y_{f}(\A)$.

\end{itemize}

\section{The graph of all pseudo-triangulations of \A}
\label{section:pseudot}


All throughout this section, $\A$ denotes a set of $n$ points in
general position in the plane,
$n_i$ of them in the  interior of $\conv(\A)$ and $n_v$  in the boundary.
Point sets in non-general position are studied in Section~\ref{section:special}.

\begin{definition} \rm
\label{defi:pseudot} A \emph{pseudo-triangle} is a simple polygon with
only three convex
vertices (called \emph{corners}) joined by three inward convex
polygonal chains (called
\emph{pseudo-edges} of the pseudo-triangle).

A \emph{pseudo-triangulation of $\A$} is a geometric non-crossing
graph  with vertex
set $\A$ and which  partitions  $\conv(\A)$  into pseudo-triangles.
\end{definition}

\begin{figure}[htb]
\begin{center} {\ }
\epsfxsize=3in\epsfbox{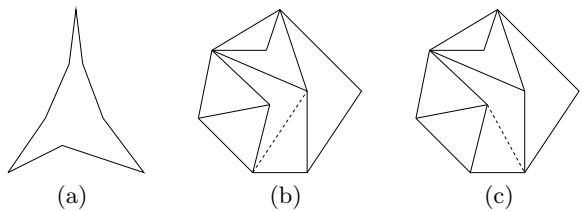} \\
{\small \hfill\qquad\quad (a) \hfill (b) \hfill (c)\qquad\quad \hfill \null\\}
\caption{(a) A pseudo-triangle. (b) A pointed pseudo-trian\-gu\-lation. (c) The
dashed edge in (b) is flipped, giving another pointed
pseudo-triangulation}\label{fig:pseudotriangs}
\end{center}
\end{figure}

Part (a) of Figure \ref{fig:pseudotriangs} shows a pseudo-triangle.
Parts (b) and (c)
show two pseudo-triangulations.

Since the maximal non-crossing graphs on $\A$ (the
\emph{triangulations} of $\A$) are
a particular case  of pseudo-triangulations, they are the maximal
pseudo-triangulations. As is well-known, they all have $2n_v+3n_i-3$
edges. It turns
out that the  pseudo-triangulations with the minimum possible number
of edges are also
very interesting from different points of view. We recall that a vertex of a
geometric graph is called \emph{pointed} if all its incident edges
span an angle
smaller than 180 degrees from that vertex. The graph itself is called
pointed if all
its vertices are pointed. The following statement comes originally from
\cite{Streinu} and a proof can also be found
in \cite{Ro-Sa-St}.

\begin{proposition}[Streinu]
\label{prop:pseudot}
Let $\A$ be a planar point set as above. Then:
\begin{enumerate}

\item[1.] Every pseudo-triangulation of $\A$ with $n_\gamma$ non-pointed
vertices and $n_\epsilon$ pointed vertices has:
$2n-3+n_\gamma =3n-3-n_\epsilon$ edges.

\item[2.] Every pointed and planar graph on $\A$ has at most $2n-3$ edges, and
is contained in some pointed pseudo-triangulation of $\A$.

\end{enumerate}
\end{proposition}

Part 1 implies that, among pseudo-triangulations of $\A$, pointed ones have the
minimum possible number of edges . For this reason they are sometimes called
\emph{minimum pseudo-triangulations}. Part 2 says that pointed
pseudo-triangulations
coincide with maximal non-crossing and pointed graphs.

Another crucial property of pseudo-triangulations is the existence of
a natural notion
of flip.  
Let $e$ be an interior edge in a pseudo-triangulation $T$ of $\A$
and let $\sigma$ be the union of the two  pseudo-triangles incident
to $e$. We regard
$\sigma$ as a graph, one of whose edges is $e$. We can consider
$\sigma\setminus e$ to
be a (perhaps degenerate) polygon, with a well-defined boundary
cycle; in degenerate
cases some edges and vertices may appear twice in the cycle. See an
example of what we
mean in Figure \ref{fig:degenerate}, in which the cycle of vertices
is $pqrstsu$ and
the cycle of edges is $pq,qr,rs,st,ts,su,up$. As in any polygon, each
(appearance of
a) vertex in the boundary cycle of $\sigma\setminus e$ is either
concave or convex. In
the figure, there are four convex vertices (corners), namely
$r$, second appearance of $s$, $u$ and $q$. Then:

\begin{figure}[htb]
\begin{center} {\ }
\epsfxsize=3in\epsfbox{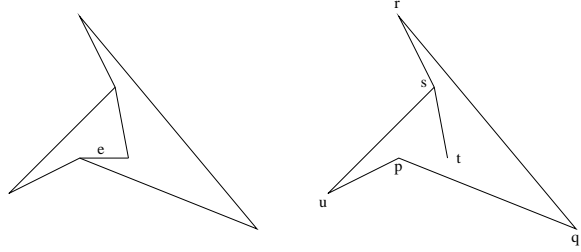}
\caption{The region $\sigma\setminus e$ is a degenerate polygon with
four corners}\label{fig:degenerate}
\end{center}
\end{figure}

\begin{lemma} \label{lemma:edge}
\begin{itemize}
\item[1.] $\sigma\setminus e$ has either 3 or 4 corners.
\item[2.] It has 3 corners if and only if exactly one of the two
end-points of $e$ is
pointed in $\sigma$. In this case $T\setminus e$ is still a
pseudo-triangulation.
\item[3.] It has 4 corners if and only if both end-points of $e$ are pointed in
$\sigma$.  In this case $T\setminus e$ is not a pseudo-triangulation and
there is a unique way to insert an edge in $T\setminus e$ to obtain another
pseudo-triangulation.
\end{itemize}
\end{lemma}

\begin{proof} Let $v_1$ and $v_2$ be the two end-points of $e$. For
each $v_i$, one of
the following three things occur: (a) $v_i$ is not-pointed in
$\sigma$, in which case
it is a corner of the two pseudo-triangles incident to $e$ and is not
a corner of
$\sigma\setminus e$; (b) $v_i$ is pointed in $\sigma$ with the big
angle exterior to
$\sigma$, in which case it is a corner of both pseudo-triangles and of
$\sigma\setminus e$ as well, or (c) $v_i$ is pointed with its big
angle interior, in
which case it is a corner in only one of the two pseudo-triangles and
not  a corner in
$\sigma\setminus e$.

In case (a), $v_i$ contributes  two more corners to the two
pseudo-triangles than to
$\sigma\setminus e$. In the other two cases, it contributes one more
corner to the
pseudo-triangles than to $\sigma\setminus e$. Since the two
pseudo-triangles have six
corners in total, $\sigma\setminus e$ has four, three or two  corners
depending on
whether both, one or none of $v_1$ and $v_2$ are pointed in $\sigma$.
The  case of two
corners is clearly impossible, which finishes the proof of 1. Part 2
only says that
``degenerate pseudo-triangles" cannot appear.

Part 3 is equivalent to saying that a pseudo-quadrangle (even a
degenerate one) can be
divided into two pseudo-triangles in exactly two ways. Indeed, these
two partitions
are obtained drawing the geodesic arcs between two opposite corners.
Such a geodesic
path consists of a unique interior edge and (perhaps) some boundary edges.
\end{proof}

Cases (2) and (3) of the above lemma will define two different types of flips
in a pseudo-triangulation.
The inverse of the first one is the insertion of an edge, in case this keeps
a pseudo-triangulation. The following statement states exactly when
this happens:

\begin{lemma} \label{lemma:insert}
Let $T$ be a pseudo-triangle with $k$ non-corners. Then, every
interior edge dividing $T$ into two pseudo-triangles makes
non-pointed exactly one
non-corner. Moreover,  there are exactly $k$ such interior edges, each making
non-pointed a different non-corner.
\end{lemma}

\begin{proof}
The first sentence follows from Lemma \ref{lemma:edge}, which says that
exactly one of the two end-points of the edge inserted is pointed (after the
insertion). For each non-corner, pointedness at the other end of the
edge implies that
the edge is the one that arises  in the geodesic arc that joins that
non-corner to the
opposite corner. This proves uniqueness and existence.
\end{proof}

\begin{figure}[htb]
\begin{center} {\ }
\epsfxsize=3in\epsfbox{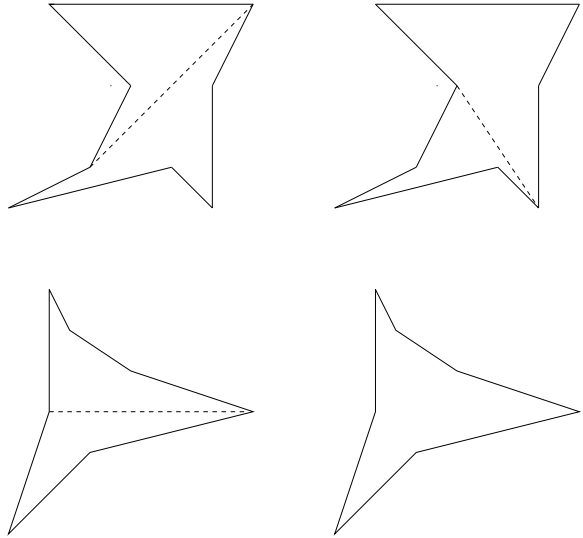}
\end{center}
\caption{Above, a diagonal-flip. Below, an insertion-deletion
flip}\label{fig:ptflips}
\end{figure}

\begin{definition}{\bf (Flips in pseudo-triangulations)} \rm
\label{defi:ptflips}
Let $T$ be a pseudo-triangu\-lation. We call \emph{flips} in $T$ the
following three
types of operations, all producing pseudo-triangulations.
See examples in Figure \ref{fig:ptflips}:

\begin{itemize}
\item[-] (\emph{Deletion flip}).  The removal of an edge $e\in T$, if
$T\setminus e$
is a pseudo-triangulation.
\item[-] (\emph{Insertion flip}).  The insertion of an edge $e\not\in
T$, if $T\cup e$
is a pseudo-triangulation.
\item[-] ({\emph{Diagonal flip})}.
   The exchange of an edge $e\in T$, if $T\setminus e$ is not a
pseudo-triangulation,
for the unique edge $e'$ such that $(T\setminus e)\cup e'$ is a
pseudo-triangulation.
\end{itemize} The \emph{graph of pseudo-triangulations} of $\A$ has
as vertices all
the pseudo-triangulations of $\A$ and as edges all flips of any of the types.
\end{definition}

\begin{proposition}\label{prop:ptgraph}
The graph of pseudo-triangulations of $\A$ is connected and regular of degree
$3n_i+n_v-3 = 3n-2n_v-3$.
\end{proposition}

\begin{proof}  There is one diagonal or deletion flip for each
interior edge, giving a
total of
$3n-3-n_\epsilon-n_v$ by Proposition \ref{prop:pseudot}. There are as
many insertion
flips as pointed interior vertices by Lemma \ref{lemma:insert},
giving $n_\epsilon
-n_v$.

To establish connectivity, let $p$ be  a point on the convex hull of $\A$. The
pseudo-triangulations of $\A$ with degree 2 at $p$
   coincide with the pseudo-triangulations of $\A\setminus\{p\}$
(together with the two
tangents from $p$ to $\A\setminus\{p\}$). By induction, we  assume all those
pseudo-triangulations to be connected in the graph.
On the other hand, in pseudo-triangulations with degree greater than
2 at $p$ all
interior edges incident to $e$ can be flipped and produce
pseudo-triangulations with
smaller degree at $e$. (Remark: if $p$ is an interior point, then  a
diagonal-flip on
an edge incident to $p$ may create another edge incident to $p$; but
for a boundary
point this cannot be the case since $p$ is a corner in the pseudo-quadrilateral
$\sigma\setminus e$ of Lemma \ref{lemma:edge}). Decreasing one by one
the number of
edges incident to $p$ will eventually lead to a pseudo-triangulation
with degree 2 at
$p$.
\end{proof}

\begin{remarks}
\rm
It is an immediate consequence of Lemma \ref{lemma:edge} that
every interior
edge in a  pointed pseudo-triangulation is flippable. This  shows
that the graph of
diagonal-flips between pointed pseudo-triangulations of $\A$ is
regular of degree
$2n_i + n_v-3$, a crucial fact in \cite{Ro-Sa-St}.

As another remark, 
one may be tempted to think that two pseudo-triangulations are
connected by a
diagonal flip  if and only if one is obtained from the other by the removal and
insertion of a single edge, but this  is  not the case: The two
pseudo-triangulations
of Figure \ref{fig:noflip}  are not connected by a diagonal flip,
according to our
definition, because the intermediate graph $T\setminus e$ is a
pseudo-triangulation.

\begin{figure}[htb]
\begin{center} {\ }
\epsfxsize=3in\epsfbox{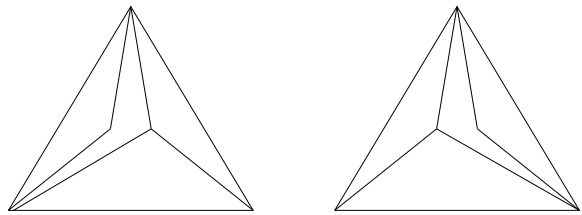}
\end{center}
\caption{These two pseudo-triangulations are not connected by a flip}
\label{fig:noflip}
\end{figure}
\end{remarks}

\subsection*{Marked non-crossing graphs on \A.}
\label{subsec:poset}

As happened with pointed pseudo-triangu\-la\-tions, Proposition \ref{prop:ptgraph}
suggests that the graph of pseudo-triangulations of $\A$ may be the
skeleton of a
simple polytope of dimension $3n_i+n_v-3$.
As a step towards this result we first look at what the face
poset of such a polytope should be. The polytope being simple means
that we want
to regard each
pseudo-triangulation $T$ as the upper bound element in a Boolean poset of order
$3n-3 - 2n_v$. This number equals, by Proposition \ref{prop:pseudot},
   the number of interior edges plus interior pointed vertices in $T$:

\begin{definition} \rm \label{defi:markedgraphs}
A \emph{marked graph} on $\A$ is a geometric graph with vertex
set $\A$ together with a subset of its vertices, that we call
``marked''.  We call a
marked graph {\em non-crossing} if it is non-crossing as a graph and
marks arise only
in pointed vertices.

We call a non-crossing marked graph
\emph{fully-marked} if it is marked at all pointed vertices. If, in
addition, it is a
pseudo-triangulation, then we call it a
\emph{fully-marked pseudo-triangulation}, abbreviated as \emph{f.m.p.t.}
\end{definition}

Marked graphs form a poset by inclusion of both the sets of edges and of marked
vertices. We say that a marked graph contains the boundary of $\A$ if
it contains all
the convex hull edges and convex hull marks. The following results
follow easily from
the corresponding statements for non-crossing graphs and
pseudo-triangulations.

\begin{proposition} \label{prop:markedpseudot} With the previous definitions:

\begin{enumerate}

\item[1.] Every marked pseudo-triangulation of $\A$ with $n_\gamma$ non-pointed
vertices, $n_\epsilon$ pointed vertices and $n_m$ marked vertices, has
   $2n-3+n_\gamma +n_m =3n-3-n_\epsilon +n_m$ edges plus marks. In
particular, all fully-marked 
pseudo-triangulations have $3n-3$ edges plus marks, $3n-3-2n_v$ of them
interior.

\item[2.]
Fully-marked pseudo-triangulations of $\A$ and maximal
non-crossing marked graphs on $\A$ are the same thing.

\item[3.] {\bf (Flips in marked pseudo-triangulations)} In a fully-marked
pseudo-triangulation of $\A$, every interior edge or interior mark
can be \emph{flipped}; once removed, there is a unique way to insert
another edge or mark to obtain a different fully-marked
pseudo-triangulation of \A.
The graph of flips between fully-marked pseudo-triangulations of
$\A$ equals the graph of pseudo-triangulations of $\A$ of Definition
\ref{defi:ptflips}.
\qed
\end{enumerate}
\end{proposition}

\begin{figure}[htb]
\begin{center} {\ }\epsfxsize=2.3in\epsfbox{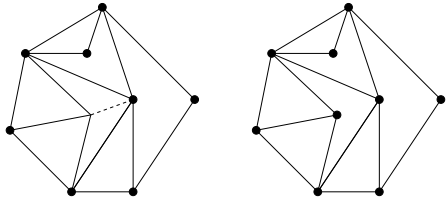}
\end{center}
\caption{Two marked pseudo-triangulations (with marks represented by
dots) related
by a flip. An edge from the left is switched to a mark on the right}
\label{fig:markedpstr}
\end{figure}

These properties imply that, if the graph of  pseudo-triangulations
of $\A$ is to be
the skeleton of a simple polytope, then the face poset of that polytope must be
(opposite to) the inclusion poset of
\emph{non-crossing marked graphs containing the boundary of  $\A$}.
Indeed, this poset
has the right ``1-skeleton'' and  the right lower ideal below every
fully-marked
pseudo-triangulation (a Boolean lattice of order $3n-3-2n_v$).



\section{The polyhedron of marked non-crossing graphs on $\A$}
\label{section:polyhedron}
In the first part of this section we do not assume $\A$ to be in 
general position.
Only after Definition \ref{defi:valid} we need general position,
among other things because we have not yet defined marked non-crossing graphs
or pseudo-triangulations for point sets in special position.
That will be done  in Section \ref{section:special}.

The setting for our construction is close to the rigid-theoretic one used in
\cite{Ro-Sa-St}.  There, the polytope to be constructed is embedded
in the space
$\R^{2n-3}$ of all infinitesimal motions of the $n$ points
$p_1,\dots,p_n$. The space
has dimension $2n-3$ because the infinitesimal motion of each point
produces two
coordinates (an infinitesimal velocity $v_i\in \R^2$) but global
translations and
rotations produce a 3-dimensional subspace of trivial motions which
are neglected.
Formally, this can be done by a quotient $\R^{2n}/ M_0$, where $M_0$ is the
3-dimensional subspace of trivial motions, or it can be done by
fixing three of the
$2n$ coordinates to be zero. For example, if the points $p_1$ and
$p_2$ do not lie in
the same horizontal line, one can take
\[ v_1^1=v_1^2=v_2^1=0.
\]
In our approach, we will consider a third coordinate $t_i$ for each
point,  related to
the ``marks'' discussed in the previous paragraphs, or to pointedness
of the vertices.


That is to say, given a set of $n$ points $\A=\{p_1,\ldots,p_n\}$ in
$\R^2$, we consider the following
$(3n-3)$-dimensional space;
\begin{equation} \label{eq:S}
\begin{matrix}
S:=\{(v_1,\ldots,v_n,t_1,\ldots,t_n)\in\left(\R^{2}\right)^n\times\R^{n}:
v_1^1=v_1^2=v_2^1=0\} \subset\R^{3n}.
\end{matrix}
\end{equation}
In it we consider the following ${n\choose 2} + n$ linear inequalities
\begin{equation}
\label{eq:v0}
H_{ij}^{+}:=\{(v,t)\in S:\langle 
p_i-p_j,v_i-v_j\rangle-|p_i-p_j|(t_i+t_j)\geq 0\}
\end{equation}
and
\begin{equation}
\label{eq:t0}
H_{0j}^{+}:=\{(v,t)\in S:t_j\geq 0\}.
\end{equation}
We denote by $H_{i,j}$ and $H_{0,j}$ their boundary hyperplanes.


\begin{definition}
\label{defi:cone}
\rm We
call \emph{expansion cone} of $\A$ and denote it $\overline{Y_0}(\A)$ 
the positive
region of the above hyperplane arrangement:
\[
\overline{Y_0}(\A):=\bigcap_{i,j\in\{0,1,\dots,n\}} H_{ij}^{+}
\]
When clear from the context we will omit the point set
$\A$ and use just $\overline{Y_0}$.
\end{definition}
Observe that the equations defining $\overline{Y_0}$ imply that for 
every $i,j$:
\[
\langle p_i-p_j, v_i-v_j \rangle \ge |p_i-p_j|(t_i+t_j) \ge 0.
\]
In particular, the vector $(v_1,\dots,v_n)$ is an expansive
infinitesimal motion of
the point set, in the standard sense.

\begin{lemma}
\label{lemma:cone}
The polyhedron $\overline{Y_0}(\A)$ has full dimension
$3n-3$ in  $S\subset \R^{3n}$ and it is a pointed polyhedral cone.
(Here, ``pointed''
means ``having the origin as a vertex'' or, equivalently,
``containing no opposite
non-zero vectors'').
\end{lemma}

\begin{proof} The vector $(v,t)$ with $v_i := p_i,
t_i:=\min_{k,l}\{|p_k-p_l|\}/4$ satisfies all the inequalities
(\ref{eq:v0}) and (\ref{eq:t0}) strictly.
In order to obtain a point in $S$ we add to it a suitable
infinitesimal trivial motion.

To prove that the cone is pointed, suppose that it contains two
opposite vectors
$(v,t)$ and $-(v,t)$. Equivalently, that $(v,t)$ lies in all the
hyperplanes $H_{i,j}$
and $H_{0,i}$. That is to say, $t_i=0$ for every $i$ and
\[
\langle v_j-v_i,p_j-p_i\rangle= 0
\]
for all $i,j$.
These last equations say that  $(v_1,\dots,v_n)$ is an infinitesimal 
flex of the complete
graph on $\A$. Since the complete graph on every full-dimensional point set
is infinitesimally rigid, $(v_1,\dots,v_n)$ is a trivial motion and equations
  (\ref{eq:S}) imply that the motion is zero.
\end{proof}

An edge $p_ip_j$ or a point $p_i$ are called {\em tight} for a certain vector
$(v,t)\in \overline{Y_0}$ if $(v,t)$ lies in the corresponding hyperplane 
$H_{i,j}$ or
$H_{0,i}$. We call {\em supporting graph} of any $(v,t)$  and denote it
$T(v,t)$ the marked graph of
tight edges for $(v,t)$ with marks at tight points for $(v,t)$.

\begin{lemma}
\label{lemma:T(v,t)}
Let $(v,t)\in\overline{Y_0}$. If $T(v,t)$ contains the boundary edges and
vertices of a
convex polygon,  then  $v_l=0$ and $t_l=0$ for every point $p_l$
in the interior of the polygon.
Therefore,  $T(v,t)$ contains the complete marked graph on the set of
vertices and interior points of the polygon.
\end{lemma}

Observe that this statement says nothing about
points in the relative interior of a boundary edge, if the polygon
has collinear points in its boundary. Indeed, such points may
have a non-zero $v_l$, namely the exterior normal to the
boundary edge containing the point.

\begin{proof}
The hypotheses are equivalent to $t_i=0$ and
$\langle p_i-p_j,v_i-v_j\rangle=0$ for all the boundary vertices
$p_i$ and boundary
edges $p_ip_j$ of the convex polygon.
We first claim that the infinitesimal expansive motion $v=(v_1,\dots,v_n)$
also preserves distances between non-consecutive
polygon vertices. Since the sum of interior angles
at vertices of an $n$-gon is independent of coordinates, a non-trivial motion
fixing the lengths of boundary edges would decrease
the interior angle at some polygon vertex $p_i$ and
then its adjacent boundary vertices get closer, what contradicts (\ref{eq:v0}).
Hence, $v$ is a translation or rotation of the polygon boundary which, by 
equations (\ref{eq:S}), is zero. On the other hand, if $v_k\neq 0$ for any $p_k$
interior to the
polygon, then $p_k$ gets closer to some boundary vertex, what using
$t_k\geq 0$ contradicts
(\ref{eq:v0}) again.

Therefore, $v_l=0$ for every point $p_l$ enclosed in the polygon (what can be
concluded from \cite[Lemma 3.2(b)]{Ro-Sa-St} as well).
Then, the equation  (\ref{eq:v0}) corresponding to $p_l$ and to  any
point $p_i$ in the
boundary of the polygon implies that $t_l\le 0$. Together with the equation
(\ref{eq:t0}) corresponding to $p_l$ this implies $t_l=0$.
\end{proof}

Obviously, $\overline{Y_0}$ is not the polyhedron
we are looking for, since its face poset does not have the desired
combinatorial structure;  it
has a unique vertex while $\A$ may have more than only one
fully-marked pseudo-triangulation. The right polyhedron for our
purposes is going
to be a convenient perturbation of $\overline{Y_0}$ obtained by translation
of its facets.

\begin{definition} \rm \label{defi:polyh}
For each $f\in \R^{n+1\choose 2}$ (with entries indexed $f_{i,j}$, for
$i,j\in\{0,\dots,n\}$) we call {\em polyhedron of expansions
constrained by $f$}, and
denote it
$\overline{Y_f}(\A)$, the polyhedron defined by the ${n\choose 2}$ equations
\begin{equation} \label{eq:vf}
\langle p_i-p_j,v_i-v_j\rangle-|p_i-p_j|(t_i+t_j)\geq f_{ij}
\end{equation}  for every $p_i,p_j \in \A$ and the $n$ equations
\begin{equation} \label{eq:tf} t_j\geq f_{0j},\qquad \forall p_i\in \A.
\end{equation}
\end{definition}

  From Lemma \ref{lemma:cone}, we conclude
that:

\begin{corollary}
\label{coro:polyhedron}
$\overline{Y_f}(\A)$ is a $(3n-3)$-dimensional unbounded polyhedron 
with at least one
vertex, for any $f$.
\qed
\end{corollary}

In the rest of this section and in Section \ref{section:valid}
we assume $\A$ to be in general position.
As before, to each feasible point $(v,t)\in \overline{Y_f}$ we associate
the marked graph
consisting of edges and vertices whose equations (\ref{eq:vf}) and
(\ref{eq:tf}) are
tight on $(v,t)$. Similarly, to a face $F$ of $\overline{Y_f}$ we associate
the tight marked
graph of any of its relative interior points. This gives an (order-reversing)
embedding of the face poset of $\overline{Y_f}$ into the poset of all
marked graphs of $\A$.
Our goal is to show that for certain choices of the constraint
parameters $f$, the
face poset of
$\overline{Y_f}$ coincides with that of non-crossing marked graphs on $\A$.

\begin{definition}
\label{defi:valid} \rm We define a choice of the constants $f$ to be
\emph{valid} if the tight marked graph
$T(F)$ of every face $F$ of $\overline{Y_f}$ is non-crossing.
\end{definition}

The proof that valid choices exist for any point set is postponed to Section
\ref{section:valid}, in order not to interrupt the current
flow of ideas. In particular, Corollary \ref{coro:exist_valid_f}
implies that the
following explicit choice is valid:

\begin{theorem} \label{thm:exist_valid_f}
The choice $f_{ij}:=
\det(O,p_i,p_j)^2$, $f_{0j}:=0$ is valid.
\end{theorem}

The main statement in the paper is then:

\begin{theorem}{\bf{(The polyhedron of marked non-crossing graphs)}}
\label{thm:polyhedron}
If $f$ is a valid choice of parameters, then
$\overline{Y_f}$ is a simple polyhedron of dimension $3n-3$
whose face poset equals (the opposite of) the poset
of non-crossing marked graphs on $\A$.
In particular:
\begin{enumerate}
\item[\textup{(a)}] Vertices of the polyhedron are in 1-to-1 correspondence
with  fully-marked pseudo-triangulations of \A.
\item[\textup{(b)}] Bounded edges correspond to flips of interior edges or
marks in fully-marked pseudo-triangulations, i.e., to fully-marked
pseudo-triangulations with one interior edge or mark removed.
\item[\textup{(c)}]
Extreme rays correspond to fully-marked
pseudo-triangulations with one convex hull edge or mark removed.
\end{enumerate}
\end{theorem}

\begin{proof}
By Corollary \ref{coro:polyhedron}, every vertex $(v,t)$ of $\overline{Y_{f}}$
has at least
$3n-3$ incident  facets. By Proposition \ref{prop:markedpseudot}, if
$f$ is valid then
the marked graph of any vertex of $\overline{Y_{f}}$ has exactly $3n-3$
edges plus marks and
is a fully-marked pseudo-triangulation. This also implies that the
polyhedron is
simple.  If we prove that all the fully-marked pseudo-triangulations appear as
vertices of $\overline{Y_{f}}$ we finish the proof, because then the face
poset of $\overline{Y_{f}}$ will have the right minimal elements and 
the right  upper
ideals of minimal
elements (the Boolean lattices of subgraphs of fully-marked
pseudo-triangulations) to
coincide with the poset of non-crossing marked graphs on $\A$.

That all fully-marked  pseudo-triangulations appear follows from
connectedness of the
graph of flips: Starting with any given vertex of $\overline{Y_{f}}$,
corresponding to a
certain f.m.p.t. $T$ of $\A$,  its $3n-3$ incident edges correspond
to the removal of
a single edge or mark in $T$.  Moreover, if the edge or mark is not
in the boundary,
Lemma \ref{lemma:T(v,t)} implies that the edge (of $\overline{Y_{f}}$)
corresponding to it is
bounded because it collapses to the origin in $\overline{Y_0}$. Then, this
edge connects the
original vertex of $\overline{Y_{f}}$ to another one which can only be the
f.m.p.t. given by
the flip in the corresponding edge or mark of $T$. Since this happens for all
vertices, and since all f.m.p.t.'s are reachable from any other one
by flips, we
conclude that they all appear as vertices.
\end{proof}

  From Theorems \ref{thm:exist_valid_f} and \ref{thm:polyhedron} it is easy
to conclude the
statements  in the introduction. The following is actually a more precise
statement implying both:

\begin{theorem}
{{\bf (The polytope of all pseudo-triangulations)}}
\label{thm:polytope}
Let $Y_f(\A)$ be the face  of $\overline{Y_f}(\A)$
defined  turning into equalities the equations (\ref{eq:vf}) and
(\ref{eq:tf}) which
correspond  to convex hull edges or convex hull points of $\A$,
and assume $f$ to be a valid choice. Then:
\begin{enumerate}
\item $Y_f(\A)$ is a simple polytope of dimension $2n_i+n-3$
whose 1-skeleton is the graph of pseudo-triangulations
of $\A$. (In particular, it is the unique maximal bounded face of 
$\overline{Y_f}(\A)$).

\item Let $F$ be the face of $Y_f(\A)$ defined by
turning into equalities the remaining equations (\ref{eq:tf}).
Then, the complement of the star of $F$ in the 
face-poset of $Y_f(\A)$ equals the poset of non-crossing graphs on $\A$
that use all the convex hull edges.
\end{enumerate}
\end{theorem}

\begin{proof}
(1) That $Y_f(\A)$ is a
bounded face follows from Lemma \ref{lemma:T(v,t)}  (it collapses to
the zero face in
$\overline{Y_0}(\A)$).
Since vertices of $\overline{Y_f}(\A)$ are f.m.p.t.'s and since all
f.m.p.t.'s contain all
the boundary edges and vertices, $Y_f(\A)$ contains all the vertices of
$\overline{Y_f}(\A)$. Hence, its vertices are in bijection with all f.m.p.t.'s
which, in turn, are
in bijection with pseudo-triangulations. Edges of
$Y_f(\A)$ correspond to f.m.p.t.'s minus one interior edge or mark, which are
precisely the flips between f.m.p.t.'s, or between pseudo-triangulations.

(2) The facets containing $F$ are those
corresponding to marks in
interior points. Then, the faces in the complement of the star of $F$
are those  in
which none of the inequalities (\ref{eq:tf}) are tight; that is to
say, they form the
poset of ``non-crossing marked graphs containing the boundary edges
and marks  but no
interior marks'', which is the same as the poset of non-crossing
graphs containing
the boundary.
\end{proof}

We now turn our attention to Theorem \ref{thm:rigidintro}. Its proof is based in the use of the 
homogeneous cone $\overline{Y_0}(\A)$ or, more preciesely, the set 
$\HH:=\{H_{ij} : i,j=1,\dots,n\}\cup \{ H_{0i} : i=1,\dots, n\}$
of hyperplanes  that define it. 

\begin{proof}[Proof of Theorem \ref{thm:rigidintro}]
Observe now that the equations defining $H_{ij}$, specialized to $t_i=0$ for every $i$, become
the equations of the infinitesimal rigidity of the complete graph on $\A$. In particular, a graph 
$G$ is rigid on $\A$ if and only if the intersection 
\[
\left( \cap_{ij\in G} H_{ij}\right)
\cap
\left( \cap_{i=1}^n H_{0i} \right)
\]
equals 0.

This happens for any pseudo-triangulation because
 Theorem \ref{thm:polyhedron} implies that
the hyperplanes corresponding to the $3n-3$
edges and marks of any fully-marked pseudo-triangulation form a basis of the (dual of)
the linear space $S$. 

To prove part (2) we only need the fact that the $3n-3$ linear hyperplanes corresponding
to a fully-marked pseudo-triangulation are independent. In particular, any subset of them is independent too. We consider the subset corresponding to the induced (marked) subgraph
on the $n-k-l$ vertices other than the $k$ pointed and $l$ non-pointed ones we are interested in.
They form an independent set  involving only $3(n-k-l)$ coordinates,
hence their number is at most $3(n-k-l)-3$ (we need to subtract 3 for the rigid motions
of the $n-k-l$ points, and here is where we need $k+l\le n-2$). 
Since the fully-marked pseudo-triangulation
has $3n-3$ edges plus marks, at least $3k+3l$ of them are incident to our subset of points.
And exactly $k$ marks are incident to our points, hence at least $2k+3l$ edges are.
\end{proof}

Actually, we can derive some consequences for general planar rigid graphs.
Observe that every planar and generically rigid graph $G$  must have between 
$2n-3$ and $3n-3$ edges (the extreme cases being an isostatic graph and a triangulation of the 
$2$-sphere). Hence, we can say that the graph $G$  has $2n-3+y$ edges, where 
and $0\le y\le n-3$. If the graph can be embedded as a pseudo-triangulation then
the embedding will have exactly $y$ non-pointed vertices. In particular, the following statement
is an indication that every planar and rigid graph can be embedded as a pseudo-triangulation:

\begin{proposition}
\label{prop.planarrigid}
Let $G$ be a planar and generically rigid graph with $n$ vertices and $2n-3+y$ edges.
Then, there is a subset $Y$ of cardinality $y$ of the vertices of $G$ such that
every set of $l$ vertices in $Y$ plus $k$ vertices not in $Y$ is incident to at least
$2k+3l$ edges, whenever $k+l\le n-2$.
\end{proposition}

\begin{proof}
Consider $G$ embedded planarly in a sufficiently generic straight-line manner. Since the embedding is planar, it can be completed to a pseudo-triangulation $T$. In particular,
the set of edges of $G$ represents an independent subset of $2n-3+y$ hyperplanes of 
$\HH$. But since the graph is rigid, adding marks to all the vertices produces a spanning set
of $3n-3+y$ hyperplanes. In between these two sets there must be a basis, consisting of the
$2n-3+y$ edges of $G$ plus $n-y$ marks. We call $Y$ the vertices not marked in this basis,
and the same argument as in the proof of Theorem \ref{thm:rigidintro} gives the statement.
\end{proof}

\begin{figure}[htb]
\begin{center} {\ }\epsfxsize=2.5in\epsfbox{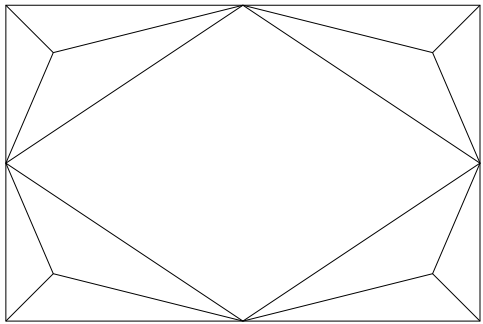}
\end{center}
\caption{A planar graph satisfying the conclusion of Proposition \ref{prop.planarrigid}
need not be rigid}
\label{fig:notrigid}
\end{figure}

It has to be said however, that a planar graph $G$ with a subset $Y$ of its vertices satisfying
Proposition \ref{prop.planarrigid} need not be generically rigid. Figure \ref{fig:notrigid} shows
an example (take as $Y$ any three of the four six-valent vertices).

\section{Valid choices of $f$}
\label{section:valid}

It remains to be proved that valid choices of parameters do exist.
In particular, that
the choice in Theorem \ref{thm:exist_valid_f} is valid. Our methods,
again inspired on
\cite{Ro-Sa-St}, give actually more: a full description of the set of
valid choices
via a set of
${n \choose 4}$ linear inequalities, one for each 4-point subset  of
the $n$ points.

\begin{definition} \rm Let $G$ be a graph embedded on $\A$, with set
of edges $E$ and
set of marked vertices $V$.  In our context, a \emph{stress} $G$ is
an assignment of
scalars
$w_{ij}$ to edges and $\alpha_j$ to marked vertices of $G$, such that for every
$(v,t)\in \R^{3n}$:
\begin{equation}\label{eq:stress}
\sum_{ij\in E} w_{ij} ( \langle p_i-p_j, v_i-v_j\rangle -|p_i-p_j|(t_i+t_j) )+
\sum_{i\in V} \alpha_i t_i= 0
\end{equation}
\end{definition}

\begin{lemma} \label{lemma:dependent} Let $\sum_{i=1}^n \lambda_i p_i=0$, $\sum
\lambda_i=0$,  be an affine dependence on a point set $\A=\{p_1,\dots,p_n\}$.
Then,
\[ w_{ij}:=\lambda_i\lambda_j \mbox{ for every $i,j$}
\] and \[
\alpha_i:=\sum_{j:ij\in E}\lambda_i\lambda_j |p_i-p_j| \mbox{ for every $i$}
\]
   defines a stress of the complete graph $G$ on \A.
\end{lemma}

\begin{proof} The condition (\ref{eq:stress}) on variables $v$ gives
\begin{equation}\label{eq:stress_in_v}
\sum_{ij\in E}w_{ij}\langle p_i-p_j, v_i-v_j\rangle = 0
\mbox{, for every } v\in \R^{2n}
\end{equation}
what can be equivalently stated as saying that the $w_{ij}$ form a
stress on the
underlying graph of $G$.
%
This is fulfilled by the $w_{ij}$'s of the
statement:
\[
\sum_{j\ne i} \lambda_i \lambda_j(p_i-p_j)=
\sum_{j=1}^n \lambda_i \lambda_j(p_i-p_j)=
\lambda_i p_i  \sum_{j=1}^n  \lambda_j - \lambda_i \sum_{j=1}^n \lambda_j p_j
= 0
\] where last equality comes from $\lambda_i$'s being an affine dependence.
Then, cancellation of the coefficient of $t_i$ in equation (\ref{eq:stress}) is
equivalent to $\alpha_i=\sum_{j:ij\in E} w_{ij} |p_i-p_j|$.
\end{proof}

Let us consider the case of four points in general position in
$\R^2$, which have a unique (up to constants) affine dependence. The
coefficients of this dependence are:
\[
\lambda_i = (-1)^{i}\det([p_1,\dots,p_4]\backslash \{p_i\})
\]
Hence, if we divide the $w_{ij}$'s and $\alpha_j$'s of the previous lemma by
the constant
\[
-\det(p_1,p_2,p_3)\det(p_1,p_2,p_4)\det(p_1,p_3,p_4)\det(p_2,p_3,p_4)
\]
   we obtain the following expressions:
\begin{equation} \label{eq:expw_ij} w_{ij}=\frac{1}{\det(p_i,p_j,p_k)
\det(p_i,p_j,p_l)},\
\alpha_{i}=\sum_{j:ij\in E} w_{ij} |p_i-p_j|
\end{equation}
where, in that of $w_{i,j}$, $k$ and $l$
denote the two indices other than $i$ and $j$.
The reason why we perform the previous rescaling is that the expressions
obtained in this way have a key property which will turn out to be
fundamental later on; see Figure \ref{fig:circuits}:

\begin{figure}[htb]
\begin{center} {\ }\epsfxsize=3in\epsfbox{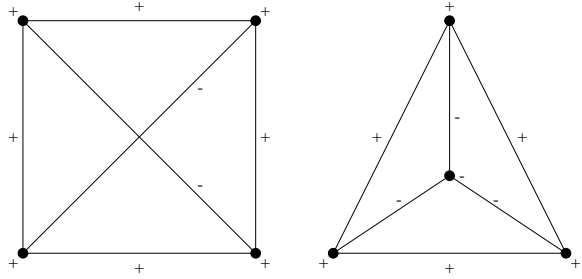}
\end{center}
\caption{The negative parts of these two marked graphs are the excluded minors
in non-crossing marked graphs of a point set in general position}
\label{fig:circuits}
\end{figure}

\begin{lemma} \label{lemma:sign_circuits} For any four
points in general position, the previous expressions give positive $w_{ij}$
and $\alpha_j$ on boundary edges and points and negative $w_{ij}$ and
$\alpha_j$ on interior edges and points.
\end{lemma}

\begin{proof} In order to check the part concerning $w_{ij}$'s we use
that $\det(q_1,q_2,q_3)$ is two times the signed area of the triangle
spanned by $q_1,q_2,q_3$: For a boundary edge the two remaining points lie on
the same side of the edge, so they have the same sign. For an interior edge,
they lie on opposite sides and therefore they have different signs.

For the $\alpha_i$'s, if $i$ is an interior point then all the
$w_{i,j}$'s in the
formula for
$\alpha_i$ are negative and, hence, $\alpha_{i}$ is also negative. If $i$ is a
boundary point then  two of the $w_{i,j}$ are positive and the third
one is negative.
But, since
\[
\sum_{j\in\{1,2,3,4\}\setminus  i} w_{i,j} (p_i-p_j) =0,
\]
the triangle inequality implies that the two positive summands
$w_{ij} |p_i-p_j|$
in the expression of $\alpha_i$ add up to a greater absolute value
than the negative
one. Hence $\alpha_i$ is positive.
\end{proof}

The previous statement is crucial to us, because no matter whether
the four points are
in convex position or one of them is inside the convex hull of the
other three, the
fully-marked pseudo-triangulations on the four points can be
characterized as the
marked graphs with nine edges plus marks and in which the missing
edge or mark is
interior (two f.m.p.t's for points in convex position, four of them
for a triangle
plus an interior point).

We conclude that:

\begin{theorem}
\label{thm:all_valids}
An $f\in \R^{n+1\choose 2}$ is valid if and only if
the following inequality holds
for every four points $\{p_1,p_2,p_3,p_4\}$ of $\A$:
\begin{equation}
\label{eq:sum_wf}
\sum_{1\leq i<j \leq 4} w_{ij} f_{ij} + \sum_{j=1}^{4} \alpha_j f_{0j}>0
\end{equation}
\label{thm:valid}
where the $w_{ij}$'s and  $\alpha_j$'s are those of (\ref{eq:expw_ij}).
\end{theorem}

\begin{proof}
Suppose first that $\A$ has only four points.
The polyhedron $\overline{Y_{f}}(\A)$ is nine-dimensional, what
implies that for every vertex $(v,t)$ of the polyhedron, the set $T(v,t)$
contains at least nine edges plus marks on those four points. Therefore,
$T(v,t)$ is the complete marked graph with an edge or mark removed.

Let us denote by $G_k$ and $G_{kl}$ the complete marked graph with a
non-marked vertex $k$ or a missing edge $kl$, respectively. Recall that by
Lemma \ref{lemma:sign_circuits} the choice of stress on four points has the
property that $G_k$ and $G_{kl}$ are fully-marked pseudo-triangulations if
and only if $\alpha_k$ and $w_{kl}$ (corresponding respectively  to the
removed mark or edge) are negative. Let us see that this is equivalent to $f$
being valid:

By the definition of stress,
\[\sum_{1\leq i < j \leq 4} w_{ij} ( \langle p_i-p_j,v_i-v_j\rangle
-|p_i-p_j|(t_i+t_j) )+
\sum_{j=1}^{4} \alpha_j t_j
\]
equals zero.
In the case of $G_k$, in which every edge and vertex except $k$ are
tight, that expression equals
\[
\sum_{1\leq i < j \leq 4} \!\! w_{ij} f_{ij} +
\sum_{j=1}^{4} \alpha_j f_{0j} +
\alpha_k (t_k-f_{0k}).
\]
In the case of $G_{kl}$, where every vertex and edge except $kl$ are tight,
it equals
\[
\sum_{1\leq i < j \leq 4} \!\! w_{ij} f_{ij} +
\sum_{j=1}^{4} \alpha_j f_{0j} +
%
w_{kl} (\langle p_k-p_l, v_k-v_l\rangle -
|p_k-p_l|(t_k+t_l) - f_{kl}).
\]

Since
$\langle p_k-p_l,v_k-v_l\rangle -|p_k-p_l|(t_k+t_l)-f_{kl} \geq 0$ and
$t_k-f_{0k}\geq 0$, by (\ref{eq:vf}) and (\ref{eq:tf}), we conclude  that in
the first and second cases above, $\alpha_k$ and $w_{kl}$ respectively are
negative if, and only if, $f$ is valid.

Now we turn to the case of a general $\A$ and our task is to prove
that a choice of
parameters $f$ is valid if and only if it is valid when restricted to
any four points.
Observe that  if $\A'\subset \A$ then $Y_f(\A')$ equals the
intersection of $Y_f(\A)$
with the subspace where $v_i=0$ and $t_i=0$ for all $p_i\in
\A\setminus \A'$. In
particular, the marked graphs on $\A'$ corresponding to faces of $Y_f(\A')$ are
subgraphs of marked graphs of faces of  $Y_f(\A)$.  Moreover,
non-crossingness of a
marked graph on $\A$ is equivalent to non-crossingness of every
induced marked graph
on four vertices: indeed, a crossing of two edges appears in the
marked graph induced
by the four end-points of the two edges,  and a non-pointed marked
vertex appears in
the marked graph induced on the  four end-points involved in any
three edges forming a
non-pointed ``letter Y'' at the non-pointed vertex.

Hence: if $f$ is valid for every four points, then none of the 4-point minors
forbidden by non-crossingness appear in faces of
$Y_f(\A)$ and $f$ is valid for $\A$. Conversely, if $f$ is not valid on some
   four point subset $\A'$, then the marked graph on $\A'$
corresponding to any vertex
of $Y_f(\A')$ would be the complete graph minus one boundary  edge or
vertex, that is
to say, it would not be non-crossing. Hence $f$ would not be valid on
$\A$ either.
\end{proof}

\begin{corollary} \label{coro:exist_valid_f}
For any $a,b\in \R^2$, the choice $f_{ij}:=
\det(a,p_i,p_j)\det(b,p_i,p_j)$, $f_{0j}:=0$ is valid.
\end{corollary}

\begin{proof}
   Consider the four points $p_i$ as fixed and regard $R:=\sum
w_{ij} f_{ij}+\sum
\alpha_j f_{0j}=\sum w_{ij} f_{ij}$ as a function of $a$ and $b$:
\[ R(a,b)=\sum_{1\leq i<j\leq 4} \det(a,p_i,p_j)\det(b,p_i,p_j) w_{ij}.
\]
We have to show that $R(a,b)$ is always positive. We actually claim
it to be always 1.
Observe first that $R(p_i,p_j)$ is trivially 1 for $i\ne j$. Since
any three of our
points are an affine basis and since $R(a,b)$ is an affine function
of $b$ for fixed
$a$, we conclude that $R(p_i,b)$ is one for every $i\in\{1,2,3,4\}$
and for every
$b$. The same argument shows that $R(a,b)$ is constantly 1: for fixed
$b$ it is an
affine function of $a$ and is equal to 1 on an affine basis.
\end{proof}

%
%

\section{Points in special position}
\label{section:special}

In this section we show that almost everything we said so far applies 
equally to point sets with collinear points. We will essentially 
follow the same steps as in Sections \ref{section:pseudot},
\ref{section:polyhedron} and \ref{section:valid}. Two subtleties are that
our definitions of pointedness or pseudo-triangulations can only be 
fully justified a posteriori,
and that the construction of the polyhedron for point sets with 
boundary collinearities is
slightly indirect: it relies in the choice of some extra exterior 
points to make colliniarities go to the interior.

\subsection*{The graph of all pseudo-triangulations of $\A$}
\label{subsec:pseudopt-sp}

\begin{definition}
\rm
\label{defi:pseudot-sp}
Let $\A$ be a finite point set in the plane, possibly with collinear points.
\begin{enumerate}
\item A graph $G$ with vertex set $\A$ is called \emph{non-crossing} 
if no edge intersects another edge of $G$
or point of $\A$ except at its end-points. In particular, if $p_1$, 
$p_2$ and $p_3$ are three
collinear points, in this order, then the edge $p_1p_3$ cannot appear 
in a non-crossing graph, independently of whether there is an edge 
incident to $p_2$ or not.

\item A \emph{pseudo-triangle} is a simple polygon with only three 
interior angles smaller than
180 degrees. A \emph{pseudo-triangulation} of $\A$ is a non-crossing 
graph with vertex set $\A$, which partitions $\conv(\A)$ into 
pseudo-triangles and such that no point in the interior of 
$\conv(\A)$ is incident to more than one angle of 180 degrees.
\end{enumerate}
\end{definition}

\begin{figure}[htb]
\begin{center}{\ }
\epsfysize=2.5in\epsfbox{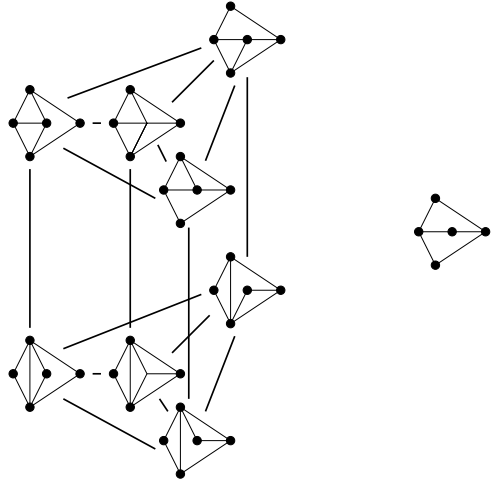}
\caption{The eight pseudo-triangulations of a point set with interior 
collinearities (left) plus
a non-crossing graph with pseudo-triangular faces but which we do not consider
a pseudo-triangulation (right)}
\label{fig:5pts}
\end{center}
\end{figure}

Figure \ref{fig:5pts} shows the eight pseudo-triangulations of a 
certain point set.
We have drawn them connected by certain flips, to be defined later, and
with certain points marked. The graph on the right of the figure is not
a pseudo-triangulation because it fails to satisfy the last condition 
in our definition. Intuitively,
the reason why we do not allow it as a pseudo-triangulation is that 
we are considering
angles of exactly 180 degrees as being
reflex, and we do not want a vertex to be incident to  two reflex angles.

\begin{figure}[htb]
\begin{center} {\ }
\epsfysize=2.5in\epsfbox{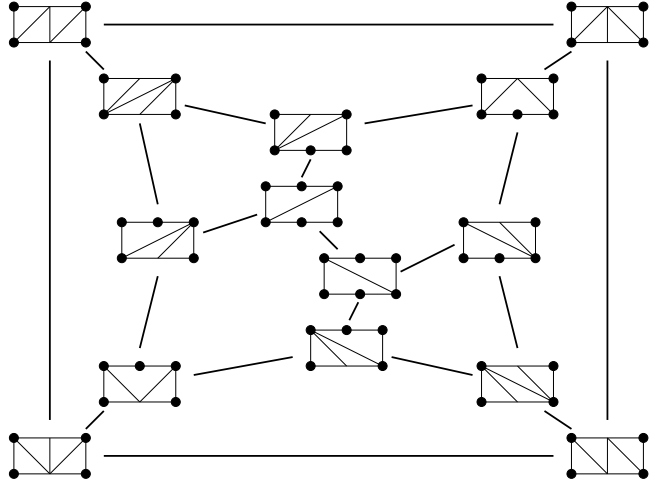}
\caption{The graph of pseudo-triangulations of a point set with 
boundary collinearities}
\label{fig:6pts}
\end{center}
\end{figure}

But if collinearities happen in the boundary of $\conv(\A)$, as in 
Figure \ref{fig:6pts}, we treat things differently. The exterior 
angle of 180 degrees is not counted as reflex, and hence
the middle point in a boundary  collinearity is allowed to be 
incident to an interior angle of
180 degrees. The following definition can be restated as
``a vertex  is pointed if and only if it is incident to a reflex 
angle'', where reflex is meant as in these last remarks.

\begin{definition}\rm
\label{defi:markedgraphs-sp}
A vertex $p$ in a non-crossing graph on $\A$ is considered
\emph{pointed} if either (1) it is a vertex of
$\conv(\A)$, (2) it is semi-interior
and not incident to any edge going through the interior of $\conv(\A)$ or
(3) it is interior  and its incident edges span at most 180 degrees.

A \emph{non-crossing marked graph} is a non-crossing graph with marks 
at some of its
pointed vertices. If all pointed vertices are marked we say the 
non-crossing graph is \emph{fully-marked}. Marks  at interior and 
semi-interior points
will be called \emph{interior marks}.
\end{definition}

For example, all the graphs of Figures \ref{fig:5pts} and 
\ref{fig:6pts} are fully-marked. That is
to say, big dots correspond exactly to pointed vertices.
Of course, fully-marked pseudo-triangulations are just 
pseudo-triangulations with marks at all their
pointed vertices.
Observe that we are calling interior marks
and edges exactly those which do not appear in all pseudo-triangulations.
 From now on, we denote by $n_i$, $n_s$ and $n_v$ the number of 
interior, semi-interior
and extremal points of $\A$.
Finally, $n=n_v+n_s+n_i$ denotes the total number of points in $\A$.
The following two statements essentially say that Proposition 
\ref{prop:markedpseudot} is valid
for non-generic configurations.

\begin{lemma}
\label{lemma:markedpseudot-sp}
Fully-marked pseudo-triangulations are exactly the maximal marked 
non-crossing graphs on $\A$.
They all have $3n-n_s-3$ edges plus marks and
$2 n_i + n -3 $ interior edges plus interior marks.
\end{lemma}

\begin{proof}
The first sentence is equivalent to saying that
  every non-crossing graph $G$ can be completed to a 
pseudo-triangulation without making any pointed vertex non-pointed. 
The proof of this is that
if $G$ is not a pseudo-triangulation then either it has a face with 
more than three corners, in which case we insert the diagonal coming 
from the geodesic between any two non-adjacent corners,
or there is an interior vertex with two angles of 180 degrees, in 
which case we choose to consider one of them as reflex and the other 
as convex, and insert the diagonal joining the convex angle to the 
opposite corner of the pseudo-triangle containing it.

To prove the cardinality of pseudo-triangulations, let $n_\epsilon$ 
denote the number of marks.
Let us think of boundary  collinearities as if they were concave 
boundary chains in our graph,
and triangulate the polygons formed by these chains by adding 
(combinatorially, or topologically)
$n_s$ edges in total. If, in addition, we consider interior angles of 
180 degrees or more as reflex and the others as convex, we get a 
graph with all the combinatorial properties of pseudo-triangulations 
and, in particular, a graph for which Proposition \ref{prop:pseudot} 
can be applied, since its
proof is purely combinatorial (a double counting of convex angles, 
combined with Euler's relation).
In particular, the extended graph has $3n-3$ edges plus marks, and 
the original graph has
  $3n-3-n_s$ of them. Since there are exactly $n_s +  n_v$ exterior 
edges and $n_v$
exterior marks in every pseudo-triangulation, the last sentence follows.
\end{proof}

\begin{lemma}
\label{lemma:ptgraph-sp}
If an interior edge or mark is removed from a fully-marked
pseudo-triangulation then there is a unique way to insert another 
edge or mark to
obtain a different fully-marked pseudo-triangulation.
\end{lemma}

\begin{proof}
If an edge is removed then there are three possibilities: (1) the 
removal does not create any new reflex angle, in which case the 
region obtained by the removal
is a pseudo-quadrangle (that is, it has four non-reflex angles), 
because the two regions incident to it had six corners in total and 
the number of them decreases by two. We insert the opposite diagonal
of it. (2) the removal creates
a new reflex angle at a vertex which was not pointed. Then the region 
obtained is a pseudo-triangle and we just add a mark at the new 
pointed vertex. (3) the removal creates a new reflex angle at a 
vertex that was already pointed. This means that after the removal 
the vertex has two reflex angles, that is to say two angles of 
exactly 180 degrees each. We insert the edge joining this vertex to 
the opposite corner of the pseudo-triangle containing the original 
reflex angle.

If a mark is removed, then the only possibility is: (4)
the pointed vertex holding the mark is incident to a unique reflex angle
(remember that we consider interior angles of 180 degrees as reflex).
We insert the edge joining the vertex to the opposite corner of the 
corresponding pseudo-triangle.
\end{proof}

\begin{definition}\rm
\label{defi:ptflip-sp}
Two fully-marked pseudo-triangulations are said to differ by a 
\emph{flip} if they differ by just one
edge or mark. Cases (1), (2), (3) and (4) in the previous proof are 
called, respectively,
  \emph{diagonal flip}, \emph{deletion flip}, \emph{mirror flip}
and \emph{insertion flip}.
\end{definition}

Of course, our definition of flips specializes to the one for points 
in general position, except that mirror flips can only appear in the 
presence of collinearities. An example of
a  mirror flip can be seen
towards the upper right corner of Figure \ref{fig:5pts}.

\begin{corollary}
\label{coro:ptgraph-sp}
The graph of flips between fully-marked pseudo-triangulations of a 
planar point set is connected and regular of degree $2n_i + n - 3$.
  \label{prop:ptgraph_deg}
\end{corollary}

The reader will have noticed that the graphs of Figures 
\ref{fig:5pts} and \ref{fig:6pts}
are more than regular of degrees 4 and 3 respectively. They are the 
graphs of  certain simple polytopes of dimensions 4 and 3. (Figure 
\ref{fig:5pts} is a prism over a simplex).

\subsection*{The case with only interior collinearities}
Now we assume that our point set $\A$ has only interior collinearities.

For each $f \in \R^{n+1}$ let $\overline{Y_f}(\A)$ be the polyhedron 
defined in Section 3.
Recall that everything we said in that section, up to Corollary 
\ref{coro:polyhedron}, is valid
for points in special position. Our main result here
is that Theorems \ref{thm:exist_valid_f} and \ref{thm:polyhedron} 
hold word by word
in the case with no boundary collinearities, except that a precision
needs to be made regarding the concept of validity.

Recall that for a given  choice of $f\in \R^{n+1\choose 2}$,
an edge $p_ip_j$ or a point $p_i$ are called \emph{tight} for
a certain $(v,t)\in \R^{3n-3}$ or for a face $F$ of $\overline{Y_f}(\A)$
  if  the corresponding equation
(\ref{eq:vf}) or (\ref{eq:tf}) is satisfied with equality.

\begin{definition}
\rm
\label{defi:valid-sp}
We call \emph{strict supporting graph} of a $(v,t)\in \R^{3n-3}$
(or  face $F$ of $\overline{Y_f}(\A)$) the marked graph of all
its tight edges and points, and denote it $T(v,t)$. We call
\emph{weak supporting graph} of a $(v,t)$ or face the marked subgraph 
consisting of
edges and points of $T(v,t)$ which define facets of $\overline{Y_f}(\A)$.

A choice of $f$ is called \emph{weakly valid} (resp., \emph{strictly 
valid}) if the
weak (resp., strict) supporting graphs of all the faces of $\overline{Y_f}(\A)$
are non-crossing marked graphs.
\end{definition}

Observe that from any weakly valid choice $f$
one can obtain strictly
valid ones: just decrease by arbitrary positive amounts the 
coordinates of $f$ corresponding to
equations which do not define facets of $\overline{Y_f}(\A)$. Hence, we could
do what follows only in terms of strict validity and would obtain the 
same polyhedron.
But weak validity is needed, as we will see in Remark 
\ref{rem:continuity}, if we want our
construction to depend continuously on the coordinates of the point set $\A$.

To obtain the equations that valid choices must satisfy we proceed as 
in Section
\ref{section:valid}. The crucial point there was that a marked graph 
is non-crossing if and only
if it does not contain the negative parts of the unique stress in 
certain subgraphs.

\begin{figure}[htb]
\begin{center} {\ }\epsfxsize=3in\epsfbox{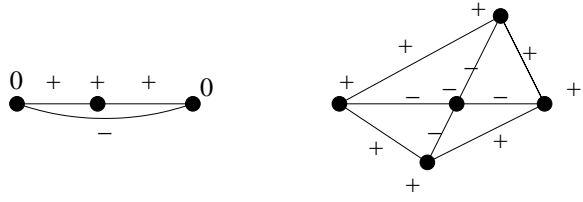}
\end{center}
\caption{The two additional excluded minors for
non-crossing marked graphs of
  a point set with interior collinearities}
\label{fig:circuits-sp}
\end{figure}

\begin{lemma}
\label{lemma:excluded}
Let $\A$ be a point set with no three collinear boundary points. 
Then, a marked graph on $\A$
is non-crossing if and only if it does not contain any of the 
following four marked subgraphs:
the negative parts of the marked
graphs displayed in Figure \ref{fig:circuits} and the negative
parts of the marked
graphs displayed in Figure \ref{fig:circuits-sp}.
\end{lemma}

\begin{proof}
Exclusion of
the negative parts of the left graphs in both figures are our 
definition of crossingness for an
unmarked graph. An interior vertex is
pointed if and only if none of the negative parts of the right graphs appear.
\end{proof}

\begin{lemma}
\label{lemma:sign_circuits-sp}
The two graphs in Figure \ref{fig:circuits-sp} have a stress with 
signs as in the figure.
\end{lemma}

\begin{proof}
For the left part it is easy to show that the following is a stress:
\[
w_{12}=\frac{1}{|p_2-p_1|}, \
w_{13}=  -\frac{1}{|p_3-p_1|}, \
w_{23}= \frac{1}{|p_3-p_2|}, \quad
\alpha_1=\alpha_3=0, \
\alpha_{2}= 2 .
\]

For the right part, observe that, by definition, stresses on a marked 
graph form a linear space.
Let the four exterior points be $p_1$,  $p_2$, $p_3$ and $p_4$, in 
cyclic order,
and let the interior point be $p_5$. We know three different stresses 
of the complete graph
on these five points: the one we used in Section \ref{section:valid} 
for the four exterior points
and the two that we have just introduced for the two collinear 
triplets. From these three
we can eliminate the coordinates of edges $p_1p_3$ and $p_2p_4$ and 
we get a stress with
the stated signs.
\end{proof}

\begin{theorem}
\label{thm:all_valids-sp}
Let $\A$ be a point set with no three collinear boundary points. 
Then, a choice of $f$
is weakly valid if it satisfies equations (\ref{eq:sum_wf}) for all 
quadruples of points
in general position plus the following sets of equations:
\begin{itemize}
\item For any  three points $p_1$, $p_2$ and $p_3$
collinear  in this order:
\begin{equation}
\label{eq:sum_wf3}
\frac{f_{12}}{|p_2-p_1|} - \frac{f_{13}}{|p_3-p_1|} + 
\frac{f_{23}}{|p_3-p_1|}  + 2 f_{02}\ge 0
\end{equation}
\item For any five points as in the right part of Figure 
\ref{fig:circuits-sp}, the following equation
where the $w_{ij}$'s and the $\alpha_i$'s form a stress with signs as 
indicated in the figure
(by convention, $w_{ij}$ equals zero for the two missing edges in the graph):
\begin{equation}
\label{eq:sum_wf5}
\sum_{1\leq i<j \leq 5} w_{ij} f_{ij} + \sum_{j=1}^{5} \alpha_j f_{0j}>0
\end{equation}
\end{itemize}
The choice is strictly valid if and only if,
moreover, the equations (\ref{eq:sum_wf3})
of all collinear triplets are satisfied strictly.
\end{theorem}

\begin{proof}
Equations (\ref{eq:sum_wf}) guarantee that no weak or strict tight 
graph contains the two excluded marked graphs of negative edges and 
points of Figure \ref{fig:circuits}. Equations (\ref{eq:sum_wf3})
and (\ref{eq:sum_wf5}) with strict inequality, do the same for the 
graphs of Figure
\ref{fig:circuits-sp}. That these equations are equivalent to strict 
validity is proved
exactly as in Section \ref{section:valid}. The reason why we allow equality in
equations (\ref{eq:sum_wf3}) if we only want a weakly valid choice is 
that the negative
part of the stress consists of a single edge. If the equation is 
satisfied with equality then
the hyperplane corresponding to this edge is a supporting hyperplane of the
face of $\overline{Y_f}$ given by the intersection of the three 
hyperplanes of the positive part
of the stress.
\end{proof}

\begin{corollary}
\label{coro:exist_valid_f-sp}
Let $\A$ be a point set with no collinear boundary points.
Any choice of $f$ satisfying  equations (\ref{eq:sum_wf}) for every 
four points in general position
plus the following ones for every collinear
triplet is weakly valid:
\begin{equation}
\label{eq:weak}
\frac{f_{12}}{|p_2-p_1|} - \frac{f_{13}}{|p_3-p_1|} + 
\frac{f_{23}}{|p_3-p_1|}  + 2 f_{02} = 0
\end{equation}

In particular,
the choices of Corollary \ref{coro:exist_valid_f} and Theorem 
\ref{thm:exist_valid_f}
are weakly valid.
\end{corollary}

\begin{proof}
For the first assertion, we need to show that equations 
(\ref{eq:sum_wf5}) follow from equations
(\ref{eq:sum_wf}) and (\ref{eq:weak}). But this is straightforward: 
from our proof
of Lemma \ref{lemma:sign_circuits-sp} it follows that
equation (\ref{eq:sum_wf5}) is just the one obtained substituting in 
(\ref{eq:sum_wf})
the values for $w_{13}$ and $w_{24}$ obtained from the two equations 
(\ref{eq:weak}).

For the last assertion, we already proved in  Corollary
\ref{coro:exist_valid_f} that the choices of $f$ introduced there 
satisfy equations
(\ref{eq:sum_wf}). It is easy, and left to the reader, to show that 
they also satisfy
  (\ref{eq:weak}).
\end{proof}

\begin{theorem}
{\bf{(Main theorem, case without boundary collinearities)}}
\label{thm:polyhedron-sp}
Let $\A$ be a point set with no three collinear points in the 
boundary of $\conv(\A)$,
and let $f$ be a weakly valid choice of parameters. Then, $\overline{Y_f}$ is a
simple polyhedron of dimension $3n-3$ with all the properties stated in
Theorems \ref{thm:polyhedron} and \ref{thm:polytope}.
\end{theorem}

\begin{proof}
Recall that if no three boundary points are collinear then every 
fully-marked pseudo-triangulation
(i.e., maximal marked non-crossing graph)
has $3n-3$ edges plus marks, exactly as in the general position case (Lemma
\ref{lemma:markedpseudot-sp}). In particular, it is still true, for 
the same reasons as in the general position case,
that $\overline{Y_f}(\A)$ is simple and all its
vertices  correspond to f.m.p.t.'s, for any valid choice of $f$. The 
rest of the arguments
in the proof of Theorem \ref{thm:polyhedron} rely on the graph of 
flips being connected, a property
that we still have. As for Theorem \ref{thm:polytope}, the face 
$Y_f(\A)$ is bounded because
Lemma \ref{lemma:T(v,t)} still applies. The rest is straightforward.
\end{proof}

\begin{remark}
\label{rem:continuity}
\rm
It is interesting to observe that taking the explicit valid choice of 
$f$ of Theorem \ref{thm:exist_valid_f}
the equations defining $\overline{Y_f}(\A)$
depend continuously on the coordinates of the points in $\A$. When three
points become collinear, the hyperplane corresponding to the (now) forbidden
edge becomes, as we said in the proof of Theorem 
\ref{thm:all_valids-sp} a supporting hyperplane
of a codimension 3 face of $\overline{Y_f}(\A)$. The combinatorics of 
the polytope changes
but maintaining its simplicity. This continuity of the defining 
hyperplanes would clearly be
impossible if we required our choice to be strictly valid for point 
sets with collinearities.
\end{remark}

\begin{example}
\rm
Let $\A_1$ and $\A_2$ be the two point sets with five points each
whose pseudo-triangulations are depicted in Figures \ref{fig:deg1} and
\ref{fig:deg2}. The first one has three collinear points and the 
second is obtained by
perturbation of the collinearity.  These two examples
were computed with the software {\tt CDD+} of Komei Fukuda 
\cite{CDD+} before we
had a clear idea of what the right definition of pseudo-triangulation 
for points in special position
should be. To emphasize the meaning of weak validity,
in Figure \ref{fig:deg1} we are showing the weak supporting graphs of 
the vertices
of $\overline{Y_f}(\A_1)$, rather than the strict ones.

\begin{figure}[htb]
\begin{center}
\epsfxsize=4.5in\epsfbox{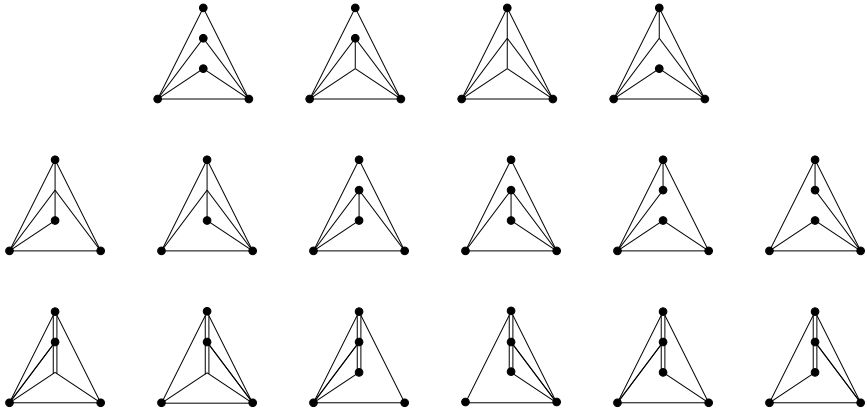}
\caption{The 16 pseudo-triangulations of $\A_1$}
\label{fig:deg1}
\end{center}
\end{figure}

\begin{figure}[htb]
\begin{center}
\epsfxsize=4.7in\epsfbox{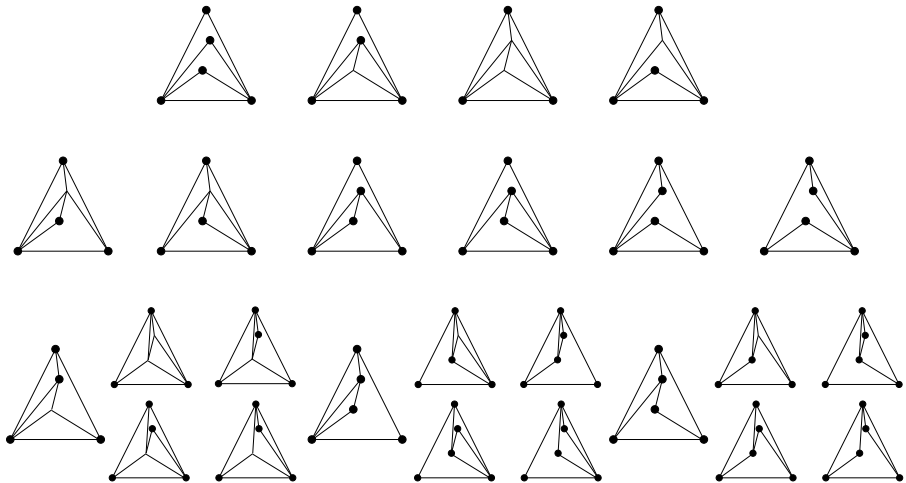}
\caption{The 25 pseudo-triangulations of $\A_2$}
\label{fig:deg2}
\end{center}
\end{figure}

The 10 first pseudo-triangulations are common to both figures (upper 
two rows). Only the
pseudo-triangulations of $\A_1$ using
the two collinear edges plus the mark at the central point of the 
collinearity are affected by the
perturbation of the point set.
This is no surprise, since these two edges plus this mark are the 
positive part of the stress
involved in the collinearity. At each of these six 
pseudo-triangulations, the hyperplane of the
big edge is tangent to the vertex of $\overline{Y_f}(\A_1)$ corresponding 
to the pseudo-triangulation.
When the collinearity is perturbed, this hyperplane moves in one of 
the two possible ways: away from the polyhedron, in which case the 
combinatorics is not changed,  or towards the interior of the 
polyhedron,
in which case the old vertex disappears and some new vertices are cut 
by this hyperplane. In our case, these two behaviors appear each in 
three of the
six ``non-strict'' pseudo-triangulations of $\A_1$. When the 
hyperplane moves towards the interior,
four new vertices appear where there was one.
\end{example}


\subsection*{Boundary collinearities}

In the presence of boundary collinearities, Lemma \ref{lemma:markedpseudot-sp}
implies a significant difference: with the equations we have used so far,
the face of $\overline{Y_0}(\A)$ defined by
tightness at boundary edges and vertices is not the origin, but an 
unbounded cone of
dimension  $n_s$. Indeed, for each semi-interior point $p_i$, the 
vector $(v,t)$
with $v_i$ an exterior
normal to the boundary of $\conv(\A)$ at $p_i$ and every other 
coordinate equal to zero
defines an extremal ray of that face. As a consequence, the 
corresponding face in
$\overline{Y_f}(\A)$ is unbounded.

We believe that it should be possible to obtain a polyhedron with the 
properties we want by just intersecting the polyhedron of our general 
definition
with $n_s$ hyperplanes. But instead of doing this
we use the following simple trick to reduce this case to the previous one.
 From a point set $\A$ with boundary collinearities
we construct another point set $\A'$ adding to $\A$ one
point in the exterior of each edge of $\conv(\A)$ that contains 
semi-interior points.

\begin{figure}[htb]
\begin{center} {\ }\epsfxsize=3in\epsfbox{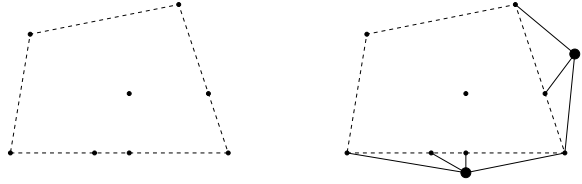}
\end{center}
\caption{The extended point set $\A'$ and the relation between 
non-crossing marked
graphs on $\A$ and $\A'$}
\label{fig:boundary}
\end{figure}

\begin{lemma}
\label{lemma:reduction}
A marked graph $G$ on $\A$ is non-crossing if and only if
it becomes a non-crossing graph on $\A'$ when we add to it
the marks on all points
of $\A'\setminus \A$ and the edges connecting each of these points to all
the points of $\A$ lying in the corresponding edge of $\conv(\A)$.
\end{lemma}

\begin{proof}
Straightforward.
\end{proof}

In particular,  we can construct the polyhedron $\overline{Y_f}(\A')$ for this
extended point set $\A'$ (taking any $f$ valid on $\A'$), and
call $\overline{Y_f}(\A)$ the face of $\overline{Y_f}(\A')$ 
corresponding to the
edges and marks mentioned in the statement of Lemma 
\ref{lemma:reduction}. Then:

\begin{corollary}
{\bf{(Main theorem, case with boundary collinearities)}}
\label{coro:reduction}
$\overline{Y_f}(\A)$
is a  simple polyhedron of dimension $3n-3-n_s$ with all the 
properties stated in
Theorem \ref{thm:polyhedron}.

Let $Y_f(\A)$ be the face of $\overline{Y_f}(\A)$
  corresponding to the $n_v+n_s$ edges between consecutive boundary points
and the $n_v$ marks at vertices of $\conv(\A)$. Let $F$ be the face of
$Y_f(\A)$ corresponding to the remaining $n-n_v$
marks. Then, $Y_f(\A)$ is a polytope of dimension $2n_i + n -3$ and $F$
is a face of it of dimension $2n_i+n_v-3$. They satisfy
  all the properties stated in Theorem  \ref{thm:polytope}.
\end{corollary}

\begin{remark}
\label{rem:semi-interior}
\rm
The reader may wonder about the combinatorics of the polyhedron
$\overline{Y_f}(\A)$ that one would obtain with the equations of the
generic case. Clearly, the tight graphs of its faces will not contain any
of the four forbidden subgraphs of Lemma \ref{lemma:excluded}. It can 
be checked
that the maximal marked graphs without those subgraphs all have 
$3n-3$ edges plus marks
and have the following characterization: as graphs they are 
pseudo-triangulations in which all the semi-interior vertices are 
incident to interior edges, and they have marks at all the boundary 
points and at the pointed interior points. In other words,
they would be the fully-marked pseudo-triangulations if we treated 
semi-interior points exactly
as interior ones, hence forbidding them to be incident to two angles 
of 180 degrees and considering
them always pointed since  they are incident to one angle of 180 degrees.

For example, in the point set of Figure \ref{fig:6pts} there are
6 such graphs, namely the ones shown in Figure \ref{fig:6pts-b}.

\begin{figure}[htb]
\begin{center} {\ }
\epsfysize=1.75in\epsfbox{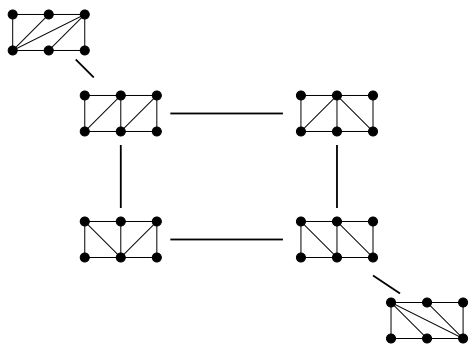}
\caption{The bounded part of $\overline{Y_f}(\A)$ for a point set with
boundary collinearities}
\label{fig:6pts-b}
\end{center}
\end{figure}

This implies that the polyhedron $\overline{Y_f}(\A)$ is still 
simple. The reason why we
prefer the definitions we have given is that the polyhedron no longer has
a unique maximal bounded face (it has three in the example of Figure 
\ref{fig:6pts-b}) and
the graph of flips is no longer regular.
\end{remark}

Observe finally that the proof of Theorem \ref{thm:rigidintro} given at the end of
Section \ref{section:polyhedron} is valid for points in special position, without much change: in all cases
the hyperplanes corresponding to the edges of a pseudo-triangulation are independent in $\overline{Y_0}(\A)$.

\end{document}